\documentclass[a4paper,12pt]{article}
\usepackage{amsmath,amsfonts,amssymb,amsthm,amscd}
                          
\renewcommand{\P}{\mathbb P}

\newcommand{\Z}{\mathbb Z}
\newcommand{\C}{\mathbb C}

\newcommand{\F}{\mathbb F}

\newcommand{\oc}{\overline c}
\newcommand{\ov}{\overline v}
\newcommand{\os}{\overline s}

\newcommand{\h}{\widehat}

\newcommand{\tD}{\widetilde D}
\newcommand{\hX}{\widehat X}
\newcommand{\hT}{\widehat T}

\newcommand{\tK}{\widetilde K}

\newcommand{\R}{\mathbb R}

\newcommand{\<}{\left <}
\renewcommand{\>}{\right >}

\renewcommand{\[}{\left [}
\renewcommand{\]}{\right ]}

\newtheorem*{theorem}{Theorem}
\newtheorem{lemma}{Lemma}
\newtheorem*{remark}{Remark}
\newenvironment{Proof}
{\noindent{\bf Proof\/}.}{{ $\square$}\smallskip\par}

\title{Existence of a singular projective variety with an arbitrary set of characteristic numbers}

\author{A. Y. Buryak \thanks{The author is partially supported by the grants RFBR-07-01-00593, NSh-709.2008.1 and the Vidi grant of NWO.}}

\date{}
\begin{document}
\maketitle
It is known that Chern characteristic numbers of compact complex manifolds cannot have arbitrary values. They satisfy certain divisability conditions.
For example (see, e.g., \cite{Hirzebruch})
\begin{eqnarray*}
  2\mid\<c_1(X),[X]\>\text{, for $dim X=1$},\\ 
  12\mid \<c_1^2(X)+c_2(X),[X]\>\text{, for $dim X=2$},\\
  24\mid \<c_1(X)c_2(X),[X]\>\text{, for $dim X=3$}.
\end{eqnarray*}

W. Ebeling and S. M. Gusein-Zade (\cite{Gusein-Zade}) offer a definition of characteristic numbers of singular compact complex analytic varieties.
For an $n$-dimensional singular analytic variety $X$, let $\nu\colon\h X\to X$ be its Nash transform and let $\hT X$ be the tautological bundle over $\h X$
(see, e.g, \cite{Gusein-Zade}).
If $X$ is embedded into a smooth complex analytic manifold $M$, then over the nonsingular part $X_{reg}$ of $X$ there is a section of $Gr_n(TM)$ given 
by the tangent space to $X$.
The Nash transform $\hX$ is the closure in $Gr_n(TM)$ of the image of this section. The bundle $\hT X$ is the restriction to $\hX$ of the tautological
bundle over $Gr_n(TM)$. 
Let the variety $X$ be compact. For a partition $I=i_1,\ldots,i_r,i_1+\ldots+i_r=n$ of $n$ the corresponding characteristic number $c_I[X]$ of the
variety $X$ is defined by
\begin{gather*}
  c_I[X]:=\<c_{i_1}(\hT X)c_{i_2}(\hT X)\cdots c_{i_r}(\hT X),\[\h X\]\>,
\end{gather*}
where $\[\h X\]$ is the fundamental class of the variety $\hX$. Let $\oc[X]$ be the vector $(c_I[X])\in\Z^{p(n)}$, where $p(n)$ is the number of 
partitions of $n$.
\begin{theorem}\label{main}
For any vector $\ov\in\Z^{p(n)}$ there exists a projective variety $X$ of 
dimension $n$
such that $\oc[X]=\ov$.
\end{theorem}
Let $V$ be an algebraic variety. R. MacPherson (\cite{MacPherson}) defined the local Euler obstruction $Eu_p(V)$ of the variety $V$ at a point $p$. 
He proved that it is a constructible function on the variety $V$. Denote this function by $Eu(X)$.
The notion of the integral with respect to the Euler characteristic was defined in \cite{Viro}. The proof of the Theorem will use the following fact.
\begin{lemma}
Let $X$ be a compact algebraic variety of dimension $n$; then $c_n[X]$ is equal to the following integral with respect to the Euler characteristic
\begin{gather*}
c_n[X]=\int\limits_X Eu(X)d\chi.
\end{gather*}
\end{lemma} 
\begin{Proof}
For any constructible function $\alpha$ on the variety $X$ R. MacPherson (\cite{MacPherson}) defined an element $c_*(\alpha)\in H_*(X)$.
From his construction it follows that 
\begin{gather*}
c_n[X]=\int\limits_X c_*(Eu(X)),
\end{gather*} 
where the integral means the degree of the class $c_*(Eu(X))$. L. Ernstr\"om (\cite{Ernstrom}) proved that for any constructible function $\alpha$ on 
a variety $X$ 
\begin{gather*}
\int\limits_X\alpha d\chi=\int\limits_X c_*(\alpha).
\end{gather*}
Lemma 1 follows from these two formulas. 
\end{Proof}
{\noindent{\bf Proof of Theorem 1\/}.}
We need some combinations of characteristic numbers (see, e.g., \cite{Milnor}).
Define two monomials in $t_1,\ldots,t_k$ to be equivalent if some permutation of 
$t_1,\ldots,t_k$ transforms 
one into the other. Define $\sum t_1^{i_1}\cdots t_r^{i_r}$ to be the sum of all 
monomials in $t_1,\cdots,t_k$ 
equivalent to $t_1^{i_1}\cdots t_r^{i_r}$. 
For any partition $I=i_1,\ldots,i_r$ of $n$, define a polynomial $s_I$ in $n$ variables as 
follows. For $k\ge n$ elementary symmetric functions $\sigma_1,\ldots,\sigma_n$ 
of $t_1,\ldots,t_k$ are algebraically independent. Let $s_I$ be the unique polynomial 
satisfying 
\begin{gather*}
  s_I(\sigma_1,\ldots,\sigma_n)=\sum t_1^{i_1}\cdots t_r^{i_r}.
\end{gather*} 
This polynomial does not depend on $k$. 
Let $F$ be a complex vector bundle over a topological space $Y$. 
For a partition $I$ of $n$ the 
cohomology class $s_I(c_1(F),\ldots,c_k(F))\in H^{2n}(Y)$ will be denoted by 
$s_I(F)$. 
For a compact analytic variety $X$ of dimension $n$ and a partition $I$ of $n$ let 
the number $s_I[X]$ be defined by
\begin{gather*}
  s_I[X]:=\<s_I(\hT X),\[\hX\]\>.
\end{gather*} 
Let $\os[X]$ be the vector $(s_I[X])\in\Z^{p(n)}$.
We have the following relationship between the vectors $\oc[X]$ and $\os[X]$ (see, e.g., \cite{Milnor}). 
There exists
a $p(n)\times p(n)$ matrix $A$ with integer coefficients and $det(A)=\pm 1$ such that, for any compact analytic 
variety $X$ of dimension $n$, one has $\oc[X]=A\os[X]$. 
Hence it is sufficient to prove that for any vector $\ov\in\Z^{p(n)}$ there exists
a projective variety $X$ such that $\os[X]=\ov$.

For two complex bundles $F,F'$ the
characteristic class $s_I(F\oplus F')$ is equal to
\begin{gather}\label{comp}
  s_I(F\oplus F')=\sum_{JK=I}s_J(F)s_K(F'),
\end{gather} 
where the sum is over all partitions $J$ and $K$ with union $JK$ equal to $I$ (\cite{Milnor}).

Let $X_1,X_2$ be two compact analytic varieties and 
$\nu_1\colon\h{X_1}\to X_1,\nu_2\colon\h{X_2}\to X_2$ be their Nash transforms. 
It is clear that the map
$(\nu_1,\nu_2)\colon\h{X_1}\times\h{X_2}\to X_1\times X_2$ is the Nash transform of 
$X_1\times X_2$. Let $p_{1,2}\colon\h{X_1}\times\h{X_2}\to\h{X_{1,2}}$ be projections; then
$\hT(X_1\times X_2)=p_1^*\hT X_1\oplus p_2^*\hT X_2$.  
Let $n_1$ and $n_2$ be the dimensions of $X_1$ and $X_2$. Let $I$ be a partition of $n_1+n_2$.
From (\ref{comp}) it follows that
\begin{gather}\label{prod}
  s_I[X_1\times X_2]=\sum_{\substack{JK=I\\|J|=n_1\\|K|=n_2}}s_J[X_1]s_K[X_2].
\end{gather} 
\begin{lemma}\label{mainlemma}
For any $i\ge 1$ there exist projective varieties $K^i_+$ and $K^i_-$ of dimension $i$
such that $s_i[K^i_{\pm}]=\pm 1$.   
\end{lemma}
We shall prove Lemma \ref{mainlemma} later.
Before that we shall deduce the statement of the Theorem from Lemma \ref{mainlemma}.
Let $J=j_1,\ldots,j_q$ be a partition of $n$. 
Let 
\begin{gather*}
  K^J_+=K^{j_1}_+\times K^{j_2}_+\times\cdots\times K^{j_q}_+,\\
  K^J_-=K^{j_1}_-\times K^{j_2}_+\times\cdots\times K^{j_q}_+.
\end{gather*}
From (\ref{prod}) it follows that 
\begin{gather*}
  s_I[K^J_{\pm}]=\sum_{\substack{I_1\cdots I_q=I\\|I_l|=j_l}}s_{I_1}[K^{j_1}_{\pm}]s_{I_2}[K^{j_2}_+]\cdots s_{I_q}[K^{j_q}_+].
\end{gather*}
A refinement of a partition $J$ means any partition which can be written as a union 
$J_1\cdots J_q$ where each $J_l$ is a partition of $j_l$. 
Consider the lexicographical order on partitions of $n$.
It is obvious that if $I$ is a refinement of $J$ then $I\le J$.
We see that the characteristic number $s_I[K^J_{\pm}]$ is zero unless the partition
$I$ is a refinement of $J$, hence
$s_I[K^J_{\pm}]=0$, if $I>J$.
We have $s_I[K^I_{\pm}]=\pm 1$.
Now it is clear that the vectors $\os[K^J_{\pm}]$ generate the whole lattice $\Z^{p(n)}$
as a semigroup. This finishes the proof of the theorem. 
{{ $\square$}\smallskip\par}
{\noindent{\bf Proof of Lemma \ref{mainlemma}\/}.}
It is known that, for any smooth compact algebraic variety $W$ of dimension $n$, there exists a smooth compact algebraic variety $V$ of dimension $n$ such
that for any partition $I$ of the number $n$ we have $c_I[V]=-c_I[W]$ (see e.g. \cite{Stong}). Denote the variety $V$ by $-W$.
We have (see e.g. \cite{Milnor}) 
\begin{gather}\label{formproj}
  s_n[\C\P^n]=n+1.  
\end{gather}
We see that existence of a variety $K^n_-$ immediately follows from existence of a variety 
$K^n_+$ because $s_n[(-\C\P^n)+nK^n_+]=-1$.
We also see that it is sufficient to construct a projective variety $\tK^n_+$ such that
$s_n[\tK^n_+]\equiv 1\mod n+1$. 

Let $n=1$. Let $\tK^1_+$ be the closure in $\C\P^2$ of the semicubical parabola $\{x^2-y^3=0\}\subset\C^2$.
From Lemma 1 and properties of the local Euler obstruction (\cite{MacPherson}) it follows that $s_1[\tK^1_+]=c_1[\tK^1_+]=3\equiv 1\mod 2$. 

Let us construct varieties $\tK^n_+$ for any $n\ge 2$.
For a smooth subvariety $X\subset\C\P^{N-1}$ of dimension $n-1$, let $CX\subset\C\P^N$ be the cone over $X$.
Let $h\in H^2(\C\P^{N-1})$ be the hyperplane class.
\begin{lemma}\label{helplemma}
Suppose the element $\left.h\right|_{X}\in H^2(X)$ is divisible by $d$; then 
\begin{gather*}
s_n[CX]\equiv ns_{n-1}[X]\mod d.
\end{gather*}
\end{lemma}
\begin{Proof}
Let $\F_{i_1,\ldots,i_s}$ be the variety consisting of flags $(V^{i_1},\ldots,V^{i_{s-1}})$ with
$V^{i_1}\subset\cdots\subset V^{i_{s-1}}\subset\C^{i_s}$ and $dim V^{i_k}=i_k$. 
Denote by $D_{i_k}$ the tautological bundle over $\F_{i_1,\ldots,i_s}$. 
Let $p$ be a point of $\C\P^N$ and let $V\subset T_p\C\P^N$ be a $d$-dimensional subspace.
Denote by $g(V)$ the unique $d$-dimensional projective subspace 
of $\C\P^N$ such that $p\in  g(V)$ and $T_p(g(V))=V$.
Let $G\subset\C\P^N$ be a $d$-dimensional projective subspace. By $k(G)$ denote the associated
$(d+1)$-dimensional vector subspace of $\C^{N+1}$.
Let $Y\subset\C\P^N$ be an $n$-dimensional subvariety. 
Consider the map
\begin{gather*} 
  \sigma\colon Y_{reg}\to\F_{1,n+1,N+1},Y_{reg}\ni p\mapsto(k(p),k(g(T_pY_{reg})))\in\F_{1,n+1,N+1}.
\end{gather*}
By definition the closure
$\overline{\sigma(Y_{reg})}$ is the Nash transform of $Y$.
The bundle $\hT Y$ is isomorphic to $\left.Hom(D_1,(D_{n+1}/D_1))\right|_{\h{Y}}$.

Let $\h{X}\subset\F_{1,n,N}$ and $\h{CX}\subset\F_{1,n+1,N+1}$ be the Nash 
transforms of $X$ and $CX$ respectively. Consider the diagram 
$$
\begin{CD}
  @.     \F_{1,2,n+1,N+1} @>\pi_2>> \F_{1,n+1,N+1}\\
  @.          @VV\pi_1V                         \\
  \F_{1,n,N} @>i>> \F_{2,n+1,N+1}    
\end{CD}
$$
where $\pi_1,\pi_2$ are the natural projections and the map $i$ is defined by                                         
\begin{gather*}
  i\colon(V^1,V^n)\mapsto(V^1\oplus k(O),V^n\oplus k(O)),
\end{gather*}
where $O\in\C\P^N$ is the vertex of the cone $CX$.
Obviously the map $i$ is injective.
Let $Y=\pi_1^{-1}(i(\h{X}))$.  
\begin{lemma}\label{lemma4}
The image of $Y$ under the map $\pi_2$ is $\h{CX}$. The map
$\pi_2\colon Y\to\h{CX}$ is birational.
\end{lemma}
\begin{Proof}
Denote by $\overline{pq}$ the line, which goes through two different points $p,q\in\C\P^N$.
From the definition of the variety $Y$ it follows that 
\begin{gather}
  Y=\{\left.(L,k(q)\oplus k(O),k(g(T_qX))\oplus k(O))\in\F_{1,2,n+1,N+1}\right|\notag\\q\in X, L\subset k(q)\oplus k(O)\}=\notag\\
  =\{\left.(k(p),k(q)\oplus k(O),k(g(T_qX))\oplus k(O))\in\F_{1,2,n+1,N+1}\right|\notag\\q\in X, p\in CX, p\in\overline{qO}\}\label{d1}.
\end{gather}
Note that, if $p\ne O$, then $q$ is uniquely determined by $p$. Denote by $Y'$ the subset of triples from (\ref{d1}) such that $p\ne O$. 

It is clear that for any point $p\in CX\backslash\{O\}$ we have 
\begin{gather*}
  k(g(T_pCX))=k(g(T_{\overline{pO}\cap X}X))\oplus k(O).  
\end{gather*} 
We see that for any element $(V^1,V^{n+1})\in\h{CX}\subset\F_{1,n+1,N+1}$ there exist points 
$p\in CX$ and $q\in X$ such that $p\in\overline{qO}$ and 
\begin{gather}\label{d2}
  V^1=k(p), V^{n+1}=k(g(T_qX))\oplus k(O).  
\end{gather}
Note that a point $q$ is not uniquely determined by the element $(V^1,V^{n+1})$.

The map $\pi_2$ just forgets the second element of the triple from (\ref{d1}) and it is clear 
that we obtain the pair $(V^1,V^{n+1})$ from (\ref{d2}). This completes the proof of the first part of the lemma.

Let $\h{CX}'=\{\left.(V^1,V^{n+1})\in\h{CX}\right|V^1\ne k(O)\}$. Note that if $(V^1,V^{n+1})\in \h{CX}'$ then a point $q$ from (\ref{d2}) is uniquely determined
and $q=\overline{pO}\cap X$. Now it is clear that the map $\pi_2$ sends $Y'$ isomorphically onto $\h{CX}'$. This concludes the proof 
of the second part of the lemma.  
\end{Proof}
By $\tD_i$ we denote the tautological bundles over $\F_{2,n+1,N+1},\F_{1,2,n+1,N+1},\F_{1,n+1,N+1}$.
By $D_i$ we denote the tautological bundles over $\F_{1,n,N}$.
We have 
\begin{gather*}
  s_n[CX]=\<s_n(\hT(CX)),\[\h{CX}\]\>=
  \<s_n(\tD^*_1\otimes(\tD_{n+1}/\tD_1)),\[\h{CX}\]\>\stackrel{\text{lemma \ref{lemma4}}}{=}\\
  =\<s_n(\tD^*_1\otimes(\tD_{n+1}/\tD_1)),[Y]\>.
\end{gather*}
The map $\pi_1\colon\F_{1,2,n+1,N+1}\to\F_{2,n+1,N+1}$ is the projectivization of the 
bundle $\tD_2$ over $\F_{2,n+1,N+1}$. We have that $i^*\tD_2\cong D_1\oplus\C$ and
$i^*\tD_{n+1}\cong D_n\oplus\C$.  
We see that the variety $Y$ is the total space $\P(D_1\oplus\C)$ of the projectivization of the bundle $D_1\oplus\C$ over $\h X$. 
By $\tau$ we denote the tautological bundle over $\P(D_1\oplus\C)$.
It is clear that $\tau=\left.\tD_1\right|_Y$. 
Therefore we have
\begin{gather*}
  \<s_n(\tD^*_1\otimes(\tD_{n+1}/\tD_1)),[Y]\>=\<s_n(\tau^*\otimes((D_n\oplus\C)/\tau)),[\P(D_1\oplus\C)]\>.
\end{gather*}
Moreover
\begin{gather*}
  s_n(\tau^*\otimes((D_n\oplus\C)/\tau))
  =s_n(\tau^*\otimes((D_n/D_1)\oplus D_1\oplus\C))=\\
  =s_n(\tau^*\otimes(D_n/D_1))+s_n(\tau^*\otimes D_1)+s_n(\tau^*)=\\
  =s_n(\tau^*\otimes D_1\otimes(D_1^*\otimes(D_n/D_1)))+s_n(\tau^*\otimes D_1)+s_n(\tau^*)=\\
  =s_n(\tau^*\otimes D_1\otimes \hT X)+s_n(\tau^*\otimes D_1)+s_n(\tau^*).
\end{gather*}
Let $c_1(\tau^*)=u\in H^2(\P(D_1\oplus\C))$.
We have $u^2=uh$.
Therefore from the assumption of the lemma it follows that for any $k\ge 2$ the element $u^k\in H^{2k}(\P(D_1\oplus\C))$ is divisible by $d$.
Hence we have 
\begin{gather*}
  \<s_n(\tau^*\otimes D_1),[\P(D_1\oplus\C)]\>=\<(u-h)^n,[\P(D_1\oplus\C)]\>\equiv 0\mod d,\\
  \<s_n(\tau^*),[\P(D_1\oplus\C)]\>=\<u^n,[\P(D_1\oplus\C)]\>\equiv 0\mod d.  
\end{gather*}
Let $x_1,\ldots,x_{n-1}$ be Chern roots of the bundle $\hT X$.
Then $x_1+u-h,\ldots,x_{n-1}+u-h$ are Chern roots of the bundle 
$\tau^*\otimes D_1\otimes \hT X$. Hence
\begin{gather*}
  \<s_n(\tau^*\otimes D_1\otimes \hT X),[\P(D_1\oplus\C)]\>=\\
  =\<\sum_{i=1}^{n-1}(x_i+u-h)^n,[\P(D_1\oplus\C)]\>\equiv\<\sum_{i=1}^{n-1}(x_i+u)^n,[\P(D_1\oplus\C)]\>\equiv\\
  \equiv\<\sum_{i=1}^{n-1}x_i^n+nu\sum_{i=1}^{n-1}x_i^{n-1},[\P(D_1\oplus\C)]\>\mod d. 
\end{gather*}
The class $\sum_{i=1}^{n-1}x_i^n\in H^{2n}(\h{X})$ is equal to zero because 
$dim_{\R}\h{X}=2n-2$. Therefore 
\begin{gather*}
  \<\sum_{i=1}^{n-1}x_i^n+nu\sum_{i=1}^{n-1}x_i^{n-1},[\P(D_1\oplus\C)]\>=
  \<nus_{n-1}(\hT X),[\P(D_1\oplus\C)]\>=\\=n\<{(\pi_{1*}}u)s_{n-1}(\hT X),\[\h{X}\]\>
  =n\<s_{n-1}(\hT X),\[\h{X}\]\>=ns_{n-1}[X].
\end{gather*}
This completes the proof of Lemma \ref{helplemma}.
\end{Proof}
Let $X=\C\P^{n-1}\hookrightarrow\C\P^{{2n\choose{n-1}}-1}$ be the image of the
Veronese embedding of degree $n+1$. Let $\tK^n_+=CX$.
From (\ref{formproj}) and lemma \ref{helplemma} it follows that $s_n[\tK^n_+]\equiv n^2\equiv 1\mod n+1$.
This concludes the proof of Lemma \ref{mainlemma}.
{{ $\square$}\smallskip\par}
\begin{remark}
It seems to be interesting to construct a cobordism theory of singular varieties associated to this notion of characteristic numbers. 
\end{remark}
{\bfseries Acknowledgements.} The author is grateful to professor S. M. Gusein-Zade for suggesting the problem and for constant attention 
to this work.

\end{document}